\documentstyle[12pt]{article}

\begin{document}
\newcommand{\ol }{\overline}
\newcommand{\ul }{\underline }
\newcommand{\ra }{\rightarrow }
\newcommand{\lra }{\longrightarrow }
\newcommand{\ga }{\gamma }
\newcommand{\st }{\stackrel }
\title{{\bf On the Nilpotent Multiplier of a Free Product}
\footnote{This research was in part supported by a grant from IPM.}}
\author{Behrooz Mashayekhy\\  Department of Mathematics, Ferdowsi University of
Mashhad,\\ P.O.Box 1159-91775, Mashhad, Iran\\ E-mail: mashaf@math.um.ac.ir\\
and\\ Institute for Studies in Theoretical Physics
and Mathematics,\\ Tehran, Iran}
\date{ }
\maketitle
\begin{abstract}
 In this paper, using a result of J. Burns and G. Ellis (Math. Z. 226(1997)
405-28.), we prove that the $c$-nilpotent multiplier (the Baer-invariant with
respect to the variety of nilpotent groups of class at most $c$, ${\cal N}_c$.)
{\it does commute} with the free product of cyclic groups of mutually coprime
order.\\
{\it A.M.S. Classification 2000}: 20E06, 20F12, 20J10\\
{\it Key words and phrases}: Nilpotent Multiplier, Baer-invariant, Free
Product\\
\end{abstract}
{\bf 1. Introduction and Motivation}

 I. Schur [12], in 1904, using projective representation theory of groups,
introduced the notion of a multiplier of a finite group. It was
known later that the Schur multiplier had a relation with homology
and cohomology of groups. In fact, if $G$ is a finite group, then
$$ M(G)\cong H^2(G,{\bf C}^{*})\ \ \ \ and\ \ \ \ \ M(G)\cong H_2(G,{\bf Z})\ ,$$
where $M(G)$ is the Schur multiplier of $G$, $H^2(G,{\bf C}^{*})$ is the second
cohomology of $G$ with coefficient in ${\bf C}^{*}$ and $H_2(G,{\bf Z})$ is the
second internal homology of $G$ [see 7].
In 1942, H. Hopf [6] proved that
$$M(G)\cong H^2(G,{\bf C}^*)\cong \frac {R\cap F'}{[R,F]}\ \ ,$$
where $G$ is presented as a quotient $G=F/R$ of a free group $F$ by a normal
subgroup $R$ in $F$. He also proved that the above formula is independent of
the presentation of $G$.

 R. Baer [1], in 1945, using the variety of groups, generalized the notion of
the Schur multiplier as follows.

 Let ${\cal V}$ be a variety of groups defined by the set of laws $V$ and let
$G$ be a group with a free presentation $ 1\lra R\lra F\lra G\lra 1$. Then the
{\it Baer-invariant} of $G$ with respect to the variety $\cal V$ is defined to
be $$ {\cal V}M(G):=\frac {R\cap V(F)}{[RV^*F]}\ , $$
where $V(F)$ is the verbal subgroup of $F$ with respect to $\cal V$ and
$$[RV^*F]=<v(f_1,\ldots ,f_{i-1},f_ir,f_{i+1},\ldots ,f_n)v(f_1,\ldots
,f_i,\ldots f_n)^{-1}\ |\ r\in R,$$ $$ 1\leq i\leq n,v\in V,f_i\in F,n\in {\bf
N}>.$$

 It is known that the Baer-invariant of a group $G$ is always abelian and
independent of the choice of the presentation of $G$. (See C. R.
Leedham-Green and S. McKay [8], from which our notation has been
taken, and H. Neumann [10] for the notion of variety of groups.)
Note that if $\cal V$ is the variety of abelian groups, $\cal A$,
then the Baer-invariant of $G$ will be
$$ {\cal A}M(G)= \frac {R\cap F'}{[R,F]}\ , $$
which is the Schur multiplier of $G$, $M(G)$. Also if ${\cal V}={\cal N}_c$ is
the variety of nilpotent groups of class at
most $c\geq 1$, then the Baer-invariant of the group $G$ with respect to ${\cal
N}_c$ will be
$$ {\cal N}_cM(G)=\frac {R\cap \ga_{c+1}(F)}{[R,\ _cF]}\ ,$$
where ${\ga}_{c+1}(F)$ is the $(c+1)$-st term of the lower
central series of $F$ and $[R, _1F]=[R,F]\ ,\ [R, _cF]=[[R,
_{c-1}F],F]$. According to J. Burns and G. Ellis' paper [2] we
shall call ${\cal N}_cM(G)$ the $c$-nilpotent multiplier of $G$
and denote it by $M^{(c)}(G)$. It is easy to see that 1-nilpotent
multiplier is actually the Schur multiplier.\\
{\bf Theorem 1.1}

 Let $\cal V$ be a variety of groups, then ${\cal V}M(-)$ is a covariant functor
from the category of all groups, ${\cal G}roups$, to the category of all abelian
groups, ${\cal  A}b$.\\
{\bf Proof.} See [8] page 107.

Now with regards to the above theorem, we are going to
concentrate on the relation between the functors, $M^{(c)}(-)$,
$c\geq 1$, and the free product as follows.

 In 1952, C. Miller [9] proved that $M(G)\cong H(G)$, where $H(G)$ is the group
of all commutator relations of $G$, taken modulo universal commutator relations.
He also showed that\\
{\bf Theorem 1.2} (C. Miller [9])

 Let $G_1$ and $G_2$ be two arbitrary groups, then $ H(G_1* G_2)\cong
H(G_1)\oplus H(G_2)$, where $G_1* G_2$ is the free product of $G_1$ and $G_2$.

 By the above theorem we can conclude the following corollary.
 \newpage
\hspace{-.2in}{\bf Corollary 1.3}

 The Schur multiplier functor, $M(-):{\cal G}roups\lra {\cal A}b$, is
coproduct-preserving. (Note that coproduct in ${\cal G}roups$ is free product
and  in ${\cal A}b$ is direct sum.)

In view of homology  and cohomology of groups, we have the following theorem.\\
{\bf Theorem 1.4}

Let $A$ be a $G$-module, then
 $H^n(-,A)\ ,\ H_n(-,A)$ are coproduct-preserving functors from ${\cal G}roups$
to ${\cal A}b$, for $n\geq 2$, i.e
$$ H^n(G_1* G_2,A)\cong H^n(G_1,A)\oplus H^n(G_2,A)\ \ \ \ \ \ for\ all\ n\geq 2\ ,$$
$$ H_n(G_1* G_2,A)\cong H_n(G_1,A)\oplus H_n(G_2,A)\ \ \ \ \ \ for\ all\ n\geq 2\ .$$
{\bf Proof.} See [5, page 220].

 Note that the above theorem does also confirm that the functor
   $$ M(-)=H_2(-,{\bf Z})=H^2(-,{\bf C}^*)\ ,$$
is coproduct-preserving.

 Now, with regards to the above theorems, it seems natural to ask whether the
$c$-nilpotent multiplier functors $M^{(c)}(-)$, $c\geq 2$, are
coproduct-preserving or not. To answer the question, first we
state an important theorem of J. Burns and G. Ellis [2,
Proposition 2.13 \&  Erratum at http://hamilton.ucg.ie/] which
is proved by a homological method.\\
{\bf Theorem 1.5} (J. Burns and G. Ellis [2])

 Let $G$ and $H$ be two arbitrary groups, then there is an isomorphism
$$M^{(2)}(G*H)\cong M^{(2)}(G)\oplus M^{(2)}(H)\oplus M(G)\otimes
H^{ab}\oplus M(H)\otimes G^{ab}\oplus Tor(G^{ab},H^{ab}) \ ,$$
where $G^{ab}=G/G'$, $H^{ab}=H/H'$ and $Tor=Tor_1^{{\bf Z}}$ .

 Now, we are ready to show that the second nilpotent multiplier functor
$M^{(2)}(-)$, {\it is not coproduct-preserving}, in general.\\
{\bf Example 1.6}

 Let $D_{\infty}=<a,b|a^2=b^2=1>\cong {\bf Z}_2*{\bf Z}_2$ be the infinite
dihedral group. Then
$$M^{(2)}(D_{\infty})\not\cong M^{(2)}({\bf Z}_2)\oplus M^{(2)}({\bf Z}_2)\
\ .$$
{\bf Proof.} By Theorem 1.5 we have
$$M^{(2)}(D_{\infty})=M^{(2)}({\bf Z}_2*{\bf Z}_2)$$ $$\cong M^{(2)}({\bf Z}_2)\oplus
M^{(2)}({\bf Z}_2) \oplus {\bf Z}_2\otimes M({\bf Z}_2)\oplus M({\bf
Z}_2)\otimes {\bf Z}_2\oplus Tor({\bf Z}_2,{\bf Z}_2) \ .$$
Clearly $M^{(2)}({\bf Z}_2)=0=M({\bf Z}_2)$ . Also it is well-known
that $ Tor({\bf Z}_2,{\bf Z}_2)\cong {\bf Z}_2\otimes {\bf Z}_2\cong {\bf Z}_2$
(see [11]). Therefore
$$ M^{(2)}({\bf Z}_2*{\bf Z}_2)\cong {\bf Z}_2 \ ,$$ but
$$ M^{(2)}({\bf Z}_2)\oplus M^{(2)}({\bf Z}_2)\cong 1\ \ .$$
Hence the result holds. $\Box $

 In spite of the above example, using Theorem 1.5, we can show that the second
nilpotent multiplier functor, $M^{(2)}(-)$, preserves the coproduct of a
finite family of cyclic groups of mutually coprime order.\\
{\bf Corollary 1.7}

 Let $\{ C_i|1\leq i\leq n\}$ be a family of cyclic groups of mutually coprime
order. Then $$M^{(2)}(\prod_{i=1}^{n}\!^{*}C_i)\cong \oplus
\sum_{i=1}^{n}M^{(2)}(C_i)\ ,$$
where $\prod_{i=1}^{n}\!^{*}C_i$ is the free product of $C_i$'s, $1\leq i\leq
n$.\\
{\bf Proof.} We proceed by induction on $n$. If $n=2$, then by Theorem 1.5 and
using the fact that the Baer-invariant of any cyclic group is trivial, we have
$$ M^{(2)}(C_1*C_2)\cong Tor(C_1,C_2)\ .$$
Since $C_1$ and $C_2$ are finite abelian groups with coprime
order, $Tor(C_1,C_2)\cong C_1\otimes C_2=1$ (see [11]).

 If $n=3$, then similarly we have
$$ M^{(2)}(C_1*C_2*C_3)\cong M^{(2)}(C_1*C_2)\oplus M^{(2)}(C_3)\oplus
M^{(1)}(C_1*C_2)\otimes C_3$$ $$\oplus (C_1*C_2)^{ab}\otimes
M^{(1)}(C_3) \oplus Tor((C_1*C_2)^{ab},C_3)$$ $$\cong
Tor(C_1\oplus C_2, C_3)\cong (C_1\oplus C_2)\otimes C_3\cong
(C_1\otimes C_3)\oplus (C_2\otimes C_3)=1\ .$$
Note that
$M^{(2)}(C_1*C_2)=M^{(2)}(C_3)=M^{(1)}(C_1*C_2)=1$
By a similar procedure we can complete the induction. $\Box $\\
{\bf 2. The Main Result}

 In this section, we are going to generalize the above corollary to the
variety of nilpotent groups of class at most $c$, ${\cal N}_c$, for all $c\geq 2$.\\
{\bf Notation 2.1}

 Let $C_i=<x_i|x_i^{r_i}>\cong {\bf Z}_{r_i}$ be cyclic group of order $r_i\ ,
\ 1\leq i\leq t$ such that $(r_i,r_j)=1$ for all $i\neq j$. Put
$C=\prod_{i=1}^{t}\!^{*}C_i$, the free product of $C_i$'s,
$1\leq i\leq t$, $F=\prod_{i=1}^{t}\!^{*}F_i$, where $F_i$ is the free group on
$\{ x_i\}$, $1\leq i\leq t$, and $S=<x_i^{r_i}|1\leq i\leq t>^F$, the normal
closure of $\{x_i^{r_i}|1\leq i\leq t\}$ in $F$. Note that $F$ is free on
$\{x_1,\ldots x_t\}$. It is easy to see that the following sequence is exact.
$$1\lra S\st{\subseteq}{\lra} F\st {nat}{\lra }C\lra 1\ \ .$$
Define by induction $\rho_1(S)=S\ ,\ \rho_{n+1}(S)=[\rho_n(S),F]$. Now by
Theorems 1.2 and 1.5, we have the following corollary.\\
{\bf Corollary 2.2}

 By the above notation and assumption, we have\\
$(i)\ \ S\cap \ga_2(F)=\rho_2(S)$.\\
$(ii)\ \ S\cap \ga_3(F)=\rho_3(S)$ and hence $\rho_2(S)\cap \ga_3(F)=\rho_3(S)$.
\\
{\bf Proof.} $(i)$ By Corollary 1.3
$M(C)=M(\prod_{i=1}^{t}\!^{*}C_i)\cong \oplus
\sum_{i=1}^{t}M(C_i)=1$. On the other hand, $M(C)\cong S\cap
\ga_2(F)/[S,F]$. Thus
$S\cap \ga_2(F)/[S,F]=1$ and so $S\cap \ga_2(F)=[S,F]=\rho_2(S)$.\\
$(ii)$ By Corollary 1.7 $M^{(2)}(C)=
M^{(2)}(\prod_{i=1}^{t}\!^{*}C_i)\cong \oplus
\sum_{i=1}^{t}M^{(2)}(C_i)=1$. Also by definition
$M^{(2)}(C)\cong S\cap \ga_3(F)/[S, _2F]$, so $\cap \ga_3(F)=[S,
_2F]=\rho_3(S)$. Moreover $\rho_3(S)\subseteq \rho_2(S)\cap
\ga_3(F)\subseteq S\cap \ga_3(F)=\rho_3(S)$ and hence
$\rho_2(S)\cap \ga_3(F)=\rho_3(S)$. $\Box$

 Now we consider the following two technical lemmas.\\
{\bf Lemma 2.3}

 By the Notation 2.1  $\rho_n(S)\cap \ga_{n+1}(F)=\rho_{n+1}(S)$, for all
$n\geq 1$.\\
{\bf Proof.} We proceed by induction on $n$. The assertion holds for $n=1,2$,
by Corollary 2.2.

 Now in order to avoid a lot of commutator manipulations, we prove the result for
$n=3$ in the special case $t=2$. Put $x=x_1\ ,\ y=x_2\ ,\ r=r_1\
,\ s=r_2$. So $F$ is free on $\{x,y\}$ and $S=<x^r,y^s>^F$.

 Let $g$ be a generator of $\rho_3(S)$, then
$$g=[(x^r)^{a_1},y^{a_2},x^{a_3}]\ or\ [(x^r)^{a_1},y^{a_2},y^{a_3}]\ or\
[(y^s)^{a_1},x^{a_2},y^{a_3}]\ or\ [(y^s),x^{a_2},x^{a_3}]\ , $$
where $a_i\in {\bf Z}$. Clearly modulo $\rho_4(S)$ we have
$$g\equiv [x^r,y,x]^{\alpha}\ or\ [x^r,y,y]^{\alpha}\ or\
[y^s,x,y]^{\alpha}\ or\ [y^s,x,x]^{\alpha}\ ,\ where\ \alpha \in {\bf Z}\ .$$

 Now, let $z\in \rho_3(S)\cap \ga_4(F)$, then $z\in \rho_3(S)$. By the above
fact and using a collecting process similar to basic commutators
(see [3]) we can obtain the following congruence modulo
$\rho_4(S)$
$$z\equiv [y^s,x,y]^{\alpha_1}[y,x^r,y]^{\beta_1}[y^s,x,x]^{\alpha_2}[y,x^r,
x]^{\beta_2}$$
$$ \equiv [y,x,y]^{s\alpha_1+r\beta_1}[y,x,x]^{s\alpha_2+r\beta_2}\ \ (mod\
\ga_4(F)),\ where\ \alpha_i,\beta_i\in {\bf Z}\ .$$
Note that we consider the order on $\{x,y\}$ as $x<y$.

 Since $z\in \rho_3(S)\cap \ga_4(F)$ and $\rho_4(S)\subseteq \ga_4(F)$, we have
$$ [y,x,y]^{s\alpha_1+r\beta_1}[y,x,x]^{s\alpha_2+r\beta_2}\in \ga_4(F)\ .$$
It is a well-known fact, by P. Hall [3, 4], that $\ga_3(F)/\ga_4(F)$ is the free
abelian group on $\{[y,x,y],[y,x,x]\}$. Therefore we conclude that
$ s\alpha_i+r\beta_i=0$, for $i=1,2$.\\
By a routine commutator calculation we have
$$ [y^s,x,y]^{\alpha_1}[y,x^r,y]^{\beta_1}\equiv
[[y^s,x]^{\alpha_1}[y,x^r]^{\beta_1},y]\ \ \ (mod\ \rho_4(S))$$
$$ [y^s,x,x]^{\alpha_2}[y,x^r,x]^{\beta_2}\equiv
[[y^s,x]^{\alpha_2}[y,x^r]^{\beta_2},x]\ \ \ (mod\ \rho_4(S)).\ \ \ (*)$$
Also
$$[y,x]^{s\alpha_i+r\beta_i}\equiv [y^s,x]^{\alpha_i}[y,x^r]^{\beta_i}\in
\rho_2(S)\ ,\ for\ i=1,2\ \ (mod\ga_3(F)).$$
since $s\alpha_i+r\beta_i=0\ ,\ i=1,2$, we have
$$ [y^s,x]^{\alpha_i}[y,x^r]^{\beta_i}\in \rho_2(S)\cap \ga_3(F)\ ,\ for\ i=1,2\
.$$
By corollary 2.2 $(ii)$ $\rho_2(S)\cap \ga_3(F)=\rho_3(S)$, thus
$$ [y^s,x]^{\alpha_i}[y,x^r]^{\beta_i}\in \rho_3(S)\ ,\ for\ i=1,2\ .$$
Therefore by $(*)$
$$ [y^s,x,y]^{\alpha_1}[y,x^r,y]^{\beta_1}\ \ ,\ \   [y^s,x,x]^{\alpha_2}[y,
x^r,x]^{\beta_2}\in \rho_4(S)).$$
Hence $z\in \rho_4(S)$, and then $\rho_3(S)\cap\ga_4(F)=\rho_4(S)$.

 Note that by a similar method we can obtain the result for $n$,
using induction hypothesis. $\Box$\\
\newpage
\hspace{-.2in}{\bf Lemma 2.4}

 By the above notation and assumption, $S\cap \ga_n(F)=\rho_n(S)$, for all $n\geq1$.\\
{\bf Proof.} We proceed by induction on $n$. For $n=1,2$ Corollary 2.2 gives the
result. Now, suppose $S\cap \ga_n(F)=\rho_(S)$ for a natural number $n$. We
show that $S\cap \ga_{n+1}(F)=\rho_{n+1}(S)$.

 Clearly $\rho_{n+1}(S)\subseteq S\cap \ga_{n+1}(F)$, also $S\cap
\ga_{n+1}(F)\subseteq S\cap \ga_n(F)=\rho_n(S)$, by induction hypothesis.
Therefore by Lemma 2.3
$$ \rho_{n+1}(S)\subseteq S\cap \ga_{n+1}(F)\subseteq \rho_n(S)\cap
\ga_{n+1}(F)=\rho_{n+1}(S)\ .$$
Hence the result holds. $\Box$

 Now, we are ready to show that the $c$-nilpotent multiplier functors, ${\cal
N}_cM(-)$, preserve the coproduct of cyclic groups of mutually coprime order,
for all $c\geq 1$.\\
{\bf Theorem 2.5}

 By the above notation and assumption,
$$ M^{(c)}(\prod_{i=1}^{t}\!^{*}C_i)\cong \oplus \sum_{i=1}^{t}M^{(c)}(C_i)=1\ ,\
for \ all\ c\geq 1\ .$$
{\bf Proof.} By Lemma 2.4 and the definition of $c$-nilpotent multiplier, we
have
$$ M^{(c)}(\prod_{i=1}^{t}\!^{*}C_i)=\frac {S\cap \ga_{c+1}(F)}{[S,
_F]}=\frac {S\cap \ga_{c+1}(F)}{\rho_{c+1}(S)}=1\ ,\ for\ all\ c\geq 1\ .$$
On the other hand, since $C_i$'s are cyclic, $M^{(c)}(C_i)=1$, so $\oplus
\sum_{i=1}^{t}{\cal N}_cM(C_i)=1$, for all $c\geq 1$. Hence the result holds.
$\Box$\\
{\bf Remark}

 In [2] it can be found some relations between the $c$-nilpotent multiplier and
the $c$-isoclinism theory of P. Hall and also the notion of $c$-capable groups.
Moreover, one may find in [2, page 423] a topological and also a homological
interpretation of the $c$-nilpotent multiplier. Thus our result, Theorem 2.5,
can be expressed and used in the above mentioned areas.

\end{document}